\documentclass[10pt]{amsart}
\usepackage{latexsym}
\usepackage{amsfonts}

\usepackage{amsmath,amsthm,amssymb}

\author[I.~Kapovich]{Ilya Kapovich}

\address{\tt Department of Mathematics, University of Illinois at
  Urbana-Champaign, 1409 West Green Street, Urbana, IL 61801, USA
  \newline http://www.math.uiuc.edu/\~{}kapovich/} \email{\tt
  kapovich@math.uiuc.edu}

\author[M.~Lustig]{Martin Lustig}\address{\tt Math\'ematiques
  (LATP), Universit\'e Paul C\'ezanne -Aix Marseille III,\\  av. Escadrille
  Normandie-Ni\'emen, 13397 Marseille 20, France} \email{\tt Martin.Lustig@univ-cezanne.fr}

\title[Domains of proper discontinuity]{Domains of proper discontinuity on the boundary of Outer space}

\newtheorem{thm}{Theorem}[section] \newtheorem{lem}[thm]{Lemma}
\newtheorem{cor}[thm]{Corollary} 
\newtheorem{prop}[thm]{Proposition} \theoremstyle{definition}
\newtheorem{defn}[thm]{Definition}
\newtheorem{notation}[thm]{Notation}
\newtheorem{conv}[thm]{Convention} \newtheorem{rem}[thm]{Remark}

\newtheorem{propdfn}[thm]{Proposition-Definition}

\def\strutdepth{\dp\strutbox}
\def \ss{\strut\vadjust{\kern-\strutdepth \sss}}
\def \sss{\vtop to \strutdepth{
\baselineskip\strutdepth\vss\llap{$\diamondsuit\;\;$}\null}}

\def\strutdepth{\dp\strutbox}
\def \sst{\strut\vadjust{\kern-\strutdepth \ssss}}
\def \ssss{\vtop to \strutdepth{
\baselineskip\strutdepth\vss\llap{$\spadesuit\;\;$}\null}}

\def\strutdepth{\dp\strutbox}
\def \ssh{\strut\vadjust{\kern-\strutdepth \sssh}}
\def \sssh{\vtop to \strutdepth{
\baselineskip\strutdepth\vss\llap{$\heartsuit\;\;$}\null}}

\newcommand{\Z}{\mathbb Z}

\def\epsilon{\varepsilon}
\def\phi{\varphi}

\newcommand{\Out}{\mbox{Out}}

\newcommand{\FN}{F_N}   % F ou F_n ou F_N ?
\newcommand{\cvn}{\mbox{cv}_N}
\newcommand{\cvnbar}{\overline{\mbox{cv}}_N}
\newcommand{\CVN}{\mbox{CV}_N}
\newcommand{\CVNbar}{\overline{\mbox{CV}}_N}
\begin{document}

\begin{abstract}
Motivated by the work of McCarthy and Papadopoulos for subgroups of
mapping class groups, we construct domains of proper discontinuity in
the compactified Outer space and in the projectivized space of
geodesic currents for any  ``sufficiently large'' subgroup of
$\Out(F_N)$ (that is, a subgroup containing a hyperbolic iwip).

As a corollary we prove that for $N\ge 3$ the action of $\Out(F_N)$ on the subset
of $\mathbb PCurr(F_N)$ consisting of all projectivized currents with
full support is properly discontinuous.
\end{abstract}

\thanks{The first author was supported by the NSF
  grant DMS-0603921}


\maketitle

\section{Introduction}

One of several important recent events in the study of mapping class
groups and Teichm\"uller spaces is the introduction and development of
the theory of convex-cocompact subgroups of mapping class groups
through the work of Farb and Mosher~\cite{FM}, Hamenstadt~\cite{Ha05},
and by Kent and Leininger~\cite{KeLe07,KeLe08,KeLe08a}. This theory is
inspired by the classical notion of convex-cocompactness for Kleinian
groups and is motivated, in part, by looking for new examples of
word-hyperbolic extensions of surface groups by non-elementary
word-hyperbolic groups.  A key component of the theory of
convex-cocompactness in the mapping class group context is the
construction of domains of discontinuity for subgroups of mapping
class groups in the Thurston boundary of the Teichm\"uller space.
This construction of domains of discontinuity was first put forward by
Masur for the handle-body group~\cite{Masur} and for arbitrary
subgroups of mapping class groups by McCarthy and
Papadopoulos~\cite{MP}, the main case being that of ``sufficiently
large" subgroups of mapping class groups. More detail about the work of
McCarthy and Papadopoulos will be given below.

For automorphisms of free groups many concepts and results proved for
mapping classes have been successfully ``translated'' in the past 20
years, sometimes giving rise to new interesting variations of the
mapping class group ideas, and occasionally even to a deeper
understanding of them. In some rare cases the innovative impulse has
even gone in the converse direction.

A satisfying translation of the subgroup theory into the $\Out(\FN)$ world, however, is still far from being achieved.  In particular, it would be much desirable if some analogues of the above quoted results for convex-cocompact subgroups could be established.  The aim of this paper is to provide a first step into this direction.  In order to state our result we first provide some terminology; more detail will be given later.

Teichm\"uller space and its Thurston boundary admit two independent translations into the $\Out(\FN)$ world: One candidate is Culler-Vogtmann's (compactified) Outer space $\CVNbar$, and the other one is the projectivized space of currents on $\FN$, which is also compact but (contrary to $\CVNbar$) infinite dimensional. The group $\Out(\FN)$ acts on both spaces, and it is known that {\em hyperbolic iwip} automorphisms (= strong analogues of pseudo-Anosov mapping classes, see Definition \ref{defn:irr} below) act on both spaces with North-South dynamics. For precise references and more details see section \ref{sec:background}.

\begin{thm}\label{thm:A}
Let $G\le \Out(F_N)$ be a subgroup
which contains at least one hyperbolic iwip. Then there exists
canonical non-empty open $G$-invariant subsets $\Delta_G^{cv} \subseteq \CVNbar$ and $\Delta_G^{curr} \subseteq \mathbb PCurr(F_N)$ on which the action of $G$ is properly discontinuous.
\end{thm}

The precise definitions of the sets $\Delta_G^{cv} \subset \CVNbar$
and $\Delta_G^{curr} \subset \mathbb PCurr(F_N)$ are given in
Section~\ref{sect:dl} and we also revisit their construction in the
discussion below.

In order to motivate and explain Theorem~\ref{thm:A} properly, we need to
recall some details of the above mentioned construction of
McCarthy and Papadopoulos~\cite{MP}.  Let $S$ be a closed hyperbolic
surface, let $Mod(S)$ be the mapping class group of $S$ and let $G\le
Mod(S)$ be a ``sufficiently large subgroup", i.e.  $G$ contains two
independent pseudo-Anosov elements. For such a $G$, there is a
well-defined \emph{limit set} $\Lambda_G\subseteq \mathcal{PML}(S)$
which is the unique smallest $G$-invariant closed subset of
$\mathcal{PML}(S)$. Here $\mathcal{PML}(S)$ is the space of projective
measured laminations on $S$. One then defines the \emph{zero locus}
$Z_G$ of $G$ as the set of all $[\lambda]\in \mathcal{PML}(S)$ such
that there exists $[\lambda']\in \Lambda_G$ satisfying
$i(\lambda,\lambda')=0$. Then $Z_G$ is a closed $G$-invariant set
containing $\Lambda_G$. Put $\Delta_G=\mathcal{PML}(S)-Z_G$. McCarthy
and Papadopoulos prove~\cite{MP} that $G$ acts properly
discontinuously on $\Delta_G$. Moreover, $G$ acts properly
discontinuously on $\mathcal T(S)\cup \Delta_G$, where $\mathcal T(S)$
is the Teichm\"uller space of $S$~\cite{KeLe08}.

In the free group context, the most frequently used analogue of
Teichm\"uller space is given by the above mentioned Outer space.  Let
$F_N$ be a free group of finite rank $N\ge 2$. The (unprojectivized)
Outer space $\cvn$ consists of all minimal free and discrete isometric
actions of $F_N$ on $\mathbb R$-trees. The projectivized Outer space
$\CVN=\mathbb P\cvn$ consists of equivalence classes of points from
$\cvn$ up to homothety. The closure $\cvnbar$ of $\cvn$ in the
equivariant Gromov-Hausdorff convergence topology consists precisely
of all the \emph{very small} minimal isometric actions of $F_N$ on
$\mathbb R$-trees, considered up to $F_N$-equivariant isometry. The
projectivization $\CVNbar=\mathbb P\cvnbar$ of $\cvnbar$ is compact
and contains $\CVN$ as a dense subset. The space $\CVNbar$ is an
analog of the Thurston compactification of the Teichm\"uller space.
The difference $\partial \CVN:=\CVNbar-\CVN$ is called the
\emph{boundary of the Outer space} $\CVN$.  All of the above spaces
come equipped with natural $\Out(F_N)$-actions, see
Section~\ref{sec:background} for more details and further references.

There is a companion space for $\cvn$ consisting of \emph{geodesic currents} on $F_N$. A geodesic current on $F_N$ is a positive Radon measure $\mu$ on $\partial^2 F_N=\partial F_N\times\partial F_N - diag$ that is invariant under the natural $F_N$-translation action and under the ``flip" map interchanging the two coordinates of $\partial^2 F_N$. Motivated by the work of Bonahon about geodesic currents on hyperbolic surfaces, geodesic currents on free groups have been introduced in a 1995 dissertation of Reiner Martin~\cite{Martin}. The notion was recently re-introduced in the work of Kapovich~\cite{Ka1,Ka2} and the theory of geodesic currents on free groups has been developed in the work of Kapovich and Lustig~\cite{KL1,KL2,KL3,KL4}, Kapovich~\cite{Ka3}, Kapovich and Nagnibeda~\cite{KN}, Francaviglia~\cite{Fra}, Coulbois, Hilion and Lustig~\cite{CHL3} and others. The space $Curr(F_N)$ of all geodesic currents is locally compact and comes equipped with a natural continuous action of $\Out(F_N)$ by linear transformations. There is a projectivization $\mathbb PCurr(F_N)$ of $Curr(F_N)$ that consists of projective classes $[\mu]$ of nonzero geodesic currents $\mu$, where two such currents are in the same projective class if they are positive scalar multiples of each other. The space $\mathbb PCurr(F_N)$ is compact and inherits a natural $\Out(F_N)$-action.
A crucial tool in this theory is the notion of a continuous \emph{geometric intersection form} $\langle\ .\ ,\ .\ \rangle: \cvnbar\times Curr(F_N)\to\mathbb R$ that was constructed by Kapovich and Lustig in~\cite{KL2}. This intersection form has some key properties in common with Bonahon's notion of a geometric intersection number between two geodesic currents on a surface, see Proposition~\ref{prop:int} below for a precise formulation.

In order to construct domains of discontinuity for subgroups of $\Out(F_N)$ in $\partial \CVN$ it turns out to be necessary to ``undualize" the picture and to play the spaces $\CVNbar$ and $\mathbb PCurr(F_N)$ off each other. We say that a subgroup $G\le \Out(F_N)$ is \emph{dynamically large} if it contains an atoroidal iwip (irreducible with irreducible powers) element, see Definition~\ref{defn:irr} for a precise definition of iwips. Being an atoroidal iwip element of $\Out(F_N)$ is the strongest free group analog of being a pseudo-Anosov mapping class. Atoroidal iwips act with ``North-South" dynamics both on $\CVNbar$ and $\mathbb PCurr(F_N)$ (see Section~\ref{sect:ns} for precise statements). This fact allows us to define in Section~\ref{sect:dl}, for a dynamically large $G\le \Out(F_N)$, its \emph{limit sets}  $\Lambda_G^{cv}\subseteq \CVNbar$ and $\Lambda_G^{curr}\subseteq \mathbb PCurr(F_N)$. When $G$ is not virtually cyclic, these limit sets are exactly the unique minimal closed $G$-invariant subsets of $\CVNbar$ and $\mathbb PCurr(F_N)$ accordingly. We then define the \emph{zero sets} of $G$:
\[
Z_G^{cv}=\{[T]\in \CVNbar : \langle T,\mu\rangle=0 \text{ for some } [\mu]\in \Lambda_G^{curr}\} \subseteq \CVNbar
\]
and
\[
Z_G^{curr}=\{[\mu]\in \mathbb PCurr(F_N) : \langle T,\mu\rangle=0 \text{ for some } [T]\in \Lambda_G^{cv}\} \subseteq \CVNbar.
\]
The zero sets are closed, $G$-invariant and contain the corresponding limit sets. We put $\Delta_G^{cv}=\CVNbar-Z_G^{cv}$ and $\Delta_G^{curr}=\mathbb PCurr(F_N)-Z_G^{curr}$, so that $\Delta_G^{cv}$ and $\Delta_G^{curr}$ are open $G$-invariant sets.

Thus Theorem~\ref{thm:A} can be viewed as a free group analogue of the
result of McCarthy and Papadopoulos~\cite{MP} for subgroups of
$Mod(S)$, mentioned above. It is easy to see that we always have
$\CVN\subseteq \Delta_G^{cv}$ because for any $[T]\in \CVN$ and any
$[\mu]\in \mathbb P Curr(F_N)$ we have $\langle T, \mu\rangle >0$. The main result of~\cite{KL3} implies a similar property for currents with full support: if $\mu\in Curr(F_N)$ is such a current then for any $[T]\in \CVNbar$ we have $\langle T, \mu\rangle >0$.  As a consequence, we obtain the following application of Theorem~\ref{thm:A}:

\begin{cor}\label{cor:B}
Let $N\ge 3$ and denote by $\mathbb PCurr_+(F_N)$ the set of all $[\mu]\in \mathbb PCurr(F_N)$ such that $\mu$ has full support. Then the action of $\Out(F_N)$ on $\mathbb PCurr_+(F_N)$ is properly discontinuous.
\end{cor}

Examples of currents with full support include the
\emph{Patterson-Sullivan currents} corresponding to points of $\cvn$,
see \cite{KN} for details.

The domains of discontinuity provided by Theorem~\ref{thm:A} are not,
in general, maximal possible. In~\cite{Gui1} Guirardel constructed a
nonempty open $\Out(F_N)$-invariant subset $\mathcal O_n\subseteq
\partial \CVN$ such that $\Out(F_N)$ acts properly discontinuously on
$\mathcal O_n$. As explained in Remark~\ref{rem:maximal} below, there
are points in $\mathcal O_n$ which do not belong to our discontinuity
domain $\Delta^{cv}_{Out(F_N)}$.

However, there is a natural class of examples where our domains of
discontinuity are likely to be maximal possible. Namely, suppose
$\phi,\psi\in \Out(F_N)$ are atoroidal iwips such that the subgroup
$\langle\phi,\psi\rangle$ is not virtually cyclic. Then, as we proved
in~\cite{KL4}, there exist $n,m\ge 1$ such that $G=\langle \phi^n,
\psi^m\rangle$ is free of rank two and such that every nontrivial
element of $G$ is again an atoroidal iwip and such that $F_N\rtimes
\langle \Phi^n,\Psi^m\rangle$ is word-hyperbolic (where $\Phi,\Psi\in
Aut(F_N)$ are representatives of $\phi$, $\psi$). Based on the mapping
class group analogy (where similar statements are known as a general
part of the theory of convex-cocompact subgroups of mapping class
groups, see~\cite{FM}), we believe that in this case it should be true
that $\Lambda_G^{cv}=Z_G^{cv}$. This would imply that
$\Delta_G^{cv}\subset\CVNbar$ is the maximal domain of discontinuity
for $G$ in this case. Moreover, again by analogy with the known mapping
class group results for convex-cocompact subgroups~\cite{FM}, we expect that in
this situation there is a natural homeomorphism between the hyperbolic
boundary of the free group $\langle \phi^n,
\psi^m\rangle$ and the limit set $\Lambda_G^{cv}$ and that, moreover,
every $[T]\in \Lambda_G^{cv}$ is uniquely ergodic, that is, there is a
unique $[\mu]\in \mathbb PCurr(F_N)$ such that $\langle
T,\mu\rangle=0$. If true, these statements would indicate that the
``Schottky type'' groups $G=\langle \phi^n,
\psi^m\rangle$ as above are good candidates for being examples of
``convex-cocompact'' subgroups of $\Out(F_N)$.

We are very grateful to Chris Leininger for explaining to us the proof
of McCarthy-Papadopoulos and for many helpful conversations.

\section{Outer space and the space of geodesic currents}\label{sec:background}

We give here only a brief overview of basic facts related to Outer space and the space of geodesic currents. We refer the reader to~\cite{CV,Ka2} for more detailed background information.

\subsection{Outer space} Let $N\ge 2$. The \emph{unprojectivized Outer
  space} $\cvn$ consists of all minimal free and discrete isometric
actions on $F_N$ on $\mathbb R$-trees (where two such actions are
considered equal if there exists an $F_N$-equivariant isometry between
the corresponding trees). There are several different topologies on
$\cvn$ that are known to coincide, in particular the equivariant
Gromov-Hausdorff convergence topology and the so-called \emph{axis} or
\emph{length function} topology. Every $T\in \cvn$ is uniquely
determined by its \emph{translation length function}
$||.||_T:F_N\to\mathbb R$, where $||g||_T$ is the translation length
of $g$ on $T$. Two trees $T_1,T_2\in \cvn$ are close if the functions
$||.||_{T_1}$ and $||.||_{T_1}$ are close pointwise on a large ball in
$F_N$. The closure $\cvnbar$ of $\cvn$ in either of these two
topologies is well-understood and known to consists precisely of all
the so-called \emph{very small} minimal isometric actions of $F_N$ on
$\mathbb R$-trees, see \cite{BF93} and \cite{CL}.  The outer
automorphism group $\Out(F_N)$ has a natural continuous right action
on $\cvnbar$ (that leaves $\cvn$ invariant) given at the level of
length functions as follows: for $T\in \cvn$ and $\phi\in \Out(F_N)$
we have $||g||_{T\phi}=||\phi(g)||_T$, where $g\in F_N$. In terms of
tree actions, $T\phi$ is equal to $T$ as a metric space, but the
action of $F_N$ is modified as: $g\underset{T\phi}{\cdot}
x=\widehat\phi(g)\underset{T}{\cdot} x$ where $x\in T$, $g\in F_N$ are
arbitrary and where $\widehat\phi\in Aut(F_N)$ is some fixed
representative of the outer automorphism $\phi$ is $Aut(F_N)$. The
\emph{projectivized Outer space} $\CVNbar=\mathbb P\cvnbar$ is defined
as the quotient $\cvn/\sim$ where for $T_1\sim T_2$ whenever
$T_2=cT_1$ for some $c>0$. One similarly defines the projectivization
$\CVNbar=\mathbb P\cvnbar$ of $\overline{cv}(F_N$ as
$\overline{cv}(F_N/\sim$ where $\sim$ is the same as above. The space
$\CVNbar$ is compact and contains $\CVN$ as a dense
$\Out(F_N)$-invariant subset. The compactification $\CVNbar$ of $\CVN$
is a free group analog of the Thurston compactification of the
Teichm\"uller space. For $T\in \cvnbar$ its $\sim$-equivalence class
is denoted by $[T]$, so that $[T]$ is the image of $T$ in $\CVNbar$.

\subsection{Geodesic currents} Let $\partial^2 F_N:=\{ (x,y)| x,y\in \partial F_N, x\ne y\}$. The action of $F_N$ by translations on its hyperbolic boundary $\partial F_N$ defines a natural diagonal action of $F_N$ on $\partial^2 F_N$. A \emph{geodesic current} on $F_N$ is a positive Radon measure on $\partial^2 F_N$ that is $F_N$-invariant and is also invariant under the ``flip" map $\partial^2 F_N\to \partial^2 F_N$, $(x,y)\mapsto (y,x)$. The space $Curr(F_N)$ of all geodesic currents on $F_N$ has a natural $\mathbb R_{\ge 0}$-linear structure and is equipped with the weak-* topology of pointwise convergence on continuous functions. Every point $T\in \cvn$ defines a \emph{simplicial chart} on $Curr(F_N)$ which allows one to think about geodesic currents as systems of nonnegative weights satisfying certain Kirchhoff-type equations; see \cite{Ka2} for details. We briefly recall the simplicial chart construction for the case where $T_A\in \cvn$ is the Cayley tree corresponding to a free basis $A$ of $F_N$. For a nondegenerate geodesic segment $\gamma=[p,q]$ in $T_A$ the \emph{two-sided cylinder} $Cyl_A(\gamma)\subseteq \partial^2 F_N$ consists of all $(x,y)\in \partial^2 F_N$ such that the geodesic from $x$ to $y$ in $T_A$ passes through $\gamma=[p,q]$. Given a nontrivial freely reduced  word $v\in F(A)=F_N$ and a current $\mu\in Curr(F_N)$, the ``weight" $\langle v,\mu\rangle_A$ is defined as $\mu(Cyl_A(\gamma))$ where $\gamma$ is any segment in the Cayley graph $T_A$ labelled by $v$ (the fact that the measure $\mu$ is $F_N$-invariant implies that a particular choice of $\gamma$ does not matter). A current $\mu$ is uniquely determined by a family of weights $\big(\langle v,\mu\rangle_A\big)_{v\in F_N-\{1\}}$. The weak-* topology on $Curr(F_N)$ corresponds to pointwise convergence of the weights for every $v\in F_N, v\ne 1$.

There is a natural left action of $\Out(F_N)$ on $Curr(F_N)$ by continuous linear transformations. Specifically, let $\mu\in Curr(F_N)$, $\phi\in \Out(F_N)$ and let $\Phi\in Aut(F_N)$ be a representative of $\phi$ in $Aut(F_N)$. Since $\Phi$ is a quasi-isometry of $F_N$, it extends to a homeomorphism of $\partial F_N$ and, diagonally, defines a homeomorphism of $\partial^2 F_N$. The measure $\phi\mu$ on $\partial^2 F_N$ is defined as follows. For a Borel subset $S\subseteq \partial^2 F_N$ we have $(\phi\mu)(S):=\mu(\Phi^{-1}(S))$. One then checks that $\phi\mu$ is a current and that it does not depend on the choice of a representative $\Phi$ of $\phi$.

For every $g\in F_N, g\ne 1$ there is an associated \emph{counting current} $\eta_g\in Curr(F_N)$. If $A$ is a free basis of $F_N$ and the conjugacy class $[g]$ of $g$ is realized by a ``cyclic word" $W$ (that is a cyclically reduced word in $F(A)$ written on a circle with no specified base-vertex), then for every nontrivial freely reduced word $v\in F(A)=F_N$ the weight $\langle v,\eta_g\rangle_A$ is equal to the total number of occurrences of $v^{\pm 1}$ in $W$ (where an occurrence of $v$ in $W$ is a vertex on $W$ such that we can read $v$ in $W$ clockwise without going off the circle). We refer the reader to \cite{Ka2} for a detailed exposition on the topic. By construction the counting current $\eta_g$ depends only on the conjugacy class $[g]$ of $[g]$ and it also satisfies $\eta_g=\eta_{g^{-1}}$. One can check~\cite{Ka2} that for $\phi\in \Out(F_N)$ and $g\in F_N, g\ne 1$ we have $\phi\eta_g=\eta_{\phi(g)}$. Scalar multiples $c\eta_g\in Curr(F_N)$, where $c\ge 0$, $g\in F_N, g\ne 1$ are called \emph{rational currents}. A key fact about $Curr(F_N)$ states that the set of all rational currents is dense in $Curr(F_N)$.

The \emph{space of projectivized geodesic currents} is defined as $\mathbb PCurr(F_N)=Curr(F_N)-\{0\}/\sim$ where $\mu_1\sim\mu_2$ whenever there exists $c>0$ such that $\mu_2=c\mu_1$. The $\sim$-equivalence class of $\mu\in Curr(F_N)-\{0\}$ is denoted by $[\mu]$.  The action of $\Out(F_N)$ on $Curr(F_N)$ descends to a continuous action of $\Out(F_N)$ on $\mathbb PCurr(F_N)$. The space $\mathbb PCurr(F_N)$ is compact and the set $\{[\eta_g\: g\in F_N, g\ne 1\}$ is a dense subset of it.

\subsection{Intersection form}

In \cite{KL2} we constructed a natural geometric \emph{intersection form} which pairs trees and currents:

\begin{prop}\label{prop:int}\cite{KL2}
Let $N\ge 2$. There exists a unique continuous map $\langle , \rangle : \cvnbar \times Curr(F_N)\to \mathbb R_{\ge 0}$ with the following properties:
\begin{enumerate}
\item We have $\langle T, c_1\mu_1+c_2\mu_2\rangle=c_1\langle T,\mu_1\rangle+c_2\langle T,\mu_2\rangle$ for any $T\in \cvnbar$, $\mu_1,\mu_2\in Curr(F_N)$, $c_1,c_2\ge 0$.
\item We have $\langle cT, \mu\rangle=c\langle T,\mu\rangle$ for any
  $T\in \cvnbar$, $\mu\in Curr(F_N)$ and $c\ge 0$.
\item We have $\langle T\phi,\mu\rangle=\langle T, \phi\mu\rangle$ for
  any $T\in \cvnbar$, $\mu\in Curr(F_N)$ and $\phi\in \Out(F_N)$.
\item We have $\langle T, \eta_g\rangle=||g||_T$ for any $T\in \cvnbar$ and $g\in F_N, g\ne 1$.
\end{enumerate}
\end{prop}

\section{North-South dynamics for atoroidal iwips}\label{sect:ns}

\begin{defn}\label{defn:irr}
An element $\phi\in \Out(F_N)$ is \emph{reducible} if there
exists a free product decomposition $F_N=C_1\ast\dots C_k\ast F'$,
where $k\ge 1$ and $C_i\ne \{1\}$, such that $\phi$ permutes the
conjugacy classes of subgroups $C_1,\dots, C_k$ in $F_N$. An element
$\phi\in \Out(F_N)$ is called \emph{irreducible} if it is not
reducible. An element $\phi\in \Out(F_N)$ is said to be
\emph{irreducible with irreducible powers} or an \emph{iwip} for
short, if for every $n\ge 1$ $\phi^n$ is irreducible (sometimes such
automorphisms are also called \emph{fully irreducible}).
\end{defn}
Thus it is easy to see that $\phi\in \Out(F_N)$ is an iwip if and only
if no positive power of $\phi$ preserves the conjugacy class of a
proper free factor of $F_N$. Recall also that $\phi\in \Out(F_N)$ is
called \emph{atoroidal} if it has no periodic conjugacy classes, that
is, if there do not exist $n\ge 1$ and $g\in F_N-\{1\}$ such that
$\phi^n$ fixes the conjugacy class $[g]$ of $g$ in $F_N$.

The following result is due to Reiner Martin~\cite{Martin}:
\begin{prop}\label{prop:rm}
Let $\phi\in \Out(F_N)$ be a atoroidal iwip. Then there exist unique $[\mu_+],[\mu_-]\in \mathbb PCurr(F_N)$ with the following properties:
\begin{enumerate}
\item The elements $[\mu_+],[\mu_-]\in \mathbb PCurr(F_N)$ are the only fixed points of $\phi$ in $\mathbb PCurr(F_N)$.
\item For any $[\mu]\ne [\mu_-]$ we have $\lim_{n\to\infty} \phi^n[\mu]=[\mu_+]$ and for any $[\mu]\ne [\mu_+]$ we have $\lim_{n\to\infty} \phi^{-n}[\mu]=[\mu_-]$.
\item We have $\phi\mu_+=\lambda_+\mu_+$ and $\phi^{-1}\mu_-=\lambda_-\mu_-$ where $\lambda_+>1$ and $\lambda_->1$. Moreover $\lambda_+$ is the Perron-Frobenius eigenvalue of any train-track representative of $\phi$ and $\lambda_-$ is the Perron-Frobenius eigenvalue of any train-track representative of $\phi^{-1}$.
\end{enumerate}
\end{prop}

A similar statement is known for $\CVNbar$ by a result of Levitt and Lustig~\cite{LL}:

\begin{prop}\label{prop:LL}
Let $\phi\in \Out(F_N)$ be an iwip. Then there exist unique $[T_+],[T_-]\in \CVNbar$ with the following properties:
\begin{enumerate}
\item The elements $[T_+],[T_-]\in \CVNbar$ are the only fixed points of $\phi$ in $\CVNbar$.
\item For any $[T]\in \CVNbar$, $[T]\ne [T_-]$ we have $\lim_{n\to\infty} [T\phi^n]=[T_+]$ and for any $[T]\in \CVNbar$, $[T]\ne [T_+]$ we have $\lim_{n\to\infty} [T\phi^{-n}]=[T_-]$.
\item We have $T_+\phi=\lambda_+T$ and $T_-\phi^{-1}=\lambda_-T_-$   where $\lambda_+>1$ and $\lambda_->1$. Moreover $\lambda_+$ is the Perron-Frobenius eigenvalue of any train-track representative of $\phi$ and $\lambda_-$ is the Perron-Frobenius eigenvalue of any train-track representative of $\phi^{-1}$.
\end{enumerate}
\end{prop}

Recall that as proved in ~\cite{KL3,KL4} we have:
\begin{prop}\label{prop:zero}
Let $\phi\in \Out(F_N)$ be an iwip and let $[T_+],[T_-]\in
\CVNbar$ and  $[\mu_+],[\mu_-]\in \mathbb PCurr(F_N)$ be as
in Proposition~\ref{prop:LL} and Proposition~\ref{prop:rm} accordingly.
Then the following hold:
\begin{enumerate}
\item For $T\in \cvnbar$ we have $\langle T,\mu_+\rangle=0$ if and only if $[T]=[T_-]$ and we have $\langle T,\mu_-\rangle=0$ if and only if $[T]=[T_+]$.
\item For $\mu\in Curr(F_N)$, $\mu\ne 0$ we have $\langle T_+,\mu\rangle=0$ if and only if $[\mu]=[\mu_-]$ and we have $\langle T_-,\mu\rangle=0$ if and only if $[\mu]=[\mu_+]$.
\item We have $\langle T_+, \mu_+\rangle>0$ and $\langle
  T_-,\mu_-\rangle>0$.
\item We have
\begin{gather*}
Stab_{\Out(F_N)}([T_+])=Stab_{\Out(F_N)}([T_-])=\\
=Stab_{\Out(F_N)}([\mu_+])=Stab_{\Out(F_N)}([\mu_-])
\end{gather*}
and this stabilizer is virtually cyclic.
\end{enumerate}
\end{prop}

We also need the following fact~\cite{KL4,KL5}:
\begin{prop}\label{prop:dl}
Let $G\le \Out(F_N)$ be a subgroup and such that there exist an atoroidal iwip $\phi\in G$. Let $[T_+(\phi)],[T_-(\phi)]\in \CVNbar$, $[\mu_+(\phi)],[\mu_-(\phi)]\in \mathbb PCurr(F_N)$ be the attracting and repelling fixed points of $\phi$ in $\CVNbar$ and $\mathbb PCurr(F_N)$ accordingly.
Then exactly one of the following occurs:
\begin{enumerate}
\item The group $G$ is virtually cyclic and preserves the sets $\{[T_+(\phi)],[T_-(\phi)]\}\subseteq \CVNbar$, $\{[\mu_+(\phi)],[\mu_-(\phi)]\}\subseteq \mathbb PCurr(F_N)$.
\item The group $G$ contains an atoroidal iwip $\psi=g\phi g^{-1}$ for
  some $g\in G$  such that
  $\{[T_+(\phi)],[T_-(\phi)]\}\cap\{[T_+(\psi)],[T_-(\psi)]\}=\emptyset$ and $\{[\mu_+(\phi)],[\mu_-(\phi)]\}\cap \{[\mu_+(\psi)],[\mu_-(\psi)]\}=\emptyset$. In this case there are some $m,n\ge 1$ such that $\langle \phi^m,\psi^n\rangle\le \Out(F_N)$ is free of rank two and, moreover, every nontrivial element of $\langle \phi^m,\psi^n\rangle$ is an atoroidal iwip.
\end{enumerate}
\end{prop}

\section{Dynamically large subgroups and their associated invariant sets}\label{sect:dl}

\begin{defn}[Dynamically large]\label{defn:dl}
We say that a subgroup $G\le \Out(F_N)$ is \emph{dynamically large} if there exist atoroidal iwip $\phi\in G$. We will say that a dynamically large subgroup $G\le \Out(F_N)$ is \emph{elementary} if it is virtually cyclic and that it is \emph{non-elementary} otherwise. Thus by Proposition~\ref{prop:dl} a non-elementary dynamically large subgroup contains a free subgroup of rank two.
\end{defn}

\begin{propdfn}[Limit set for a nonelementary dynamically large subgroup]
Let $G\le \Out(F_N)$ be a non-elementary dynamically large subgroup. Then the following hold:
\begin{enumerate}
\item There exists a unique minimal nonempty closed $G$-invariant subset $\Lambda_G^{cv}$ of $\CVNbar$. Moreover, $\Lambda_G^{cv}\subseteq \partial \CVN$ and  every $G$-orbit of a point of $\Lambda_G^{cv}$ is dense in $\Lambda_G^{cv}$. We call $\Lambda_G^{cv}$ the \emph{limit set} of $G$ in $\CVNbar$.

\item There exists a unique minimal nonempty closed $G$-invariant subset $\Lambda_G^{curr}$ of $\mathbb PCurr(F_N)$. Moreover, every $G$-orbit of a point of $\Lambda_G^{curr}$ is dense in $\Lambda_G^{curr}$. We call $\Lambda_G^{curr}$ the \emph{limit set} of $G$ in $\mathbb PCurr(F_N)$.
\end{enumerate}
\end{propdfn}
\begin{proof}

We will prove part (1) of the proposition as the proof of (2) is similar.

By Proposition~\ref{prop:dl} and Proposition~\ref{prop:zero} there exist atoriodal iwips $\phi,\psi\in G$, such that $[T_\pm(\phi)],[T_\pm(\psi)]$ are four distinct points and $[\mu_\pm(\phi)], [\mu_\pm(\psi)]$ are four distinct points.

Put $\Lambda_G^{cv}$ to be the closure in $\CVNbar$ of the orbit $G[T_+(\phi)]$. Thus $\Lambda_G^{cv}$ is a closed $G$-invariant subset.
Let $X\subseteq \CVNbar$ be a nonempty closed invariant subset. Note that $X$ must contain a point $[T]\ne [T_\pm(\phi)]$ since $G$ is non-elementary and $\psi\in G$ does not leave invariant a nonempty subset of $\{[T_\pm(\phi)]\}$. Then $\lim_{n\to\infty} \phi^n[T]=[T_+(\phi)]$ and hence $[T_+\phi]\in X$. Since $X$ is closed and $G$-invariant, it follows that  $\Lambda_G^{cv}=\overline{G[T_+(\phi)]}\subseteq X$. Thus $\Lambda_G^{cv}=\overline{G[T_+(\phi)]}$ is the unique minimal nonempty closed $G$-invariant subset of $\CVNbar$.
It follows that the $G$-orbit of every point of $\Lambda_G^{cv}$ is dense in $\Lambda_G^{cv}$, since the closure of such an orbit is a closed $G$-invariant set and thus must contain $\Lambda_G^{cv}$.
\end{proof}

\begin{defn}[Limit set of an elementary dynamically large subgroup]
Let $G\le \Out(F_N)$ be an elementary dynamically large subgroup and let $\phi\in G$ be an atoroidal iwip. We put $\Lambda_G^{cv}:=\{[T_+(\phi)], [T_-(\phi)]\}$ and call it the \emph{limit set} of $G$ in $\CVNbar$. Similarly, we put $\Lambda_G^{curr}:=\{[\mu_+(\phi)], [\mu_-(\phi)]\}$ and call it the \emph{limit set} of $G$ in $\mathbb PCurr(F_N)$.
\end{defn}

\begin{defn}[Zero sets]
Let $G\le \Out(F_N)$ be a dynamically large subgroup.

\begin{enumerate}
\item Put
\[
Z_G^{cv}=\{[T]\in \CVNbar: \langle T,\mu\rangle=0 \text{ for some } [\mu]\in \Lambda_G^{curr}\}.
\]
\item Put
\[
Z_G^{curr}=\{[\mu]\in \mathbb PCurr(F_N): \langle T,\mu\rangle=0 \text{ for some } [T]\in \Lambda_G^{cv}\}.
\]
\end{enumerate}
\end{defn}

The following is an immediate corollary of the definitions and of Proposition~\ref{prop:zero}:

\begin{prop}\label{prop:z-elementary}
Let $G\le \Out(F_N)$ be an elementary dynamically large subgroup and let $\phi\in G$ be an atoroidal iwip. Then $Z_G^{cv}=\Lambda_G^{cv}=\{[T_\pm (\phi)]\}$ and $Z_G^{curr}=\Lambda_G^{curr}=\{[\mu_\pm (\phi)]\}$.
\end{prop}

The proof of the following proposition is straightforward, particularly in view of continuity of the intersection form, and we leave the details to the reader.
\begin{prop}\label{prop:Z}
Let $G\le \Out(F_N)$. Then $Z_G^{cv}\subseteq \CVNbar$ and $Z_G^{curr}\subseteq \mathbb PCurr(F_N)$ are closed $G$-invariant subsets and $\Lambda_G^{cv}\subseteq \Z_G^{cv}$ and $\Lambda_G^{curr}\subseteq Z_G^{curr}$.
\end{prop}

\begin{rem}\label{rem:poles}
Let $G\le \Out(F_N)$ and let $\phi\in G$ be an atoroidal iwip. Then $[T_{\pm }(\phi)]\in \Lambda_G^{cv}$ and $[\mu_\pm (\phi)]\in \Lambda_G^{curr}$. Indeed, if $G$ is elementary, this follows from the definitions. If $G$ is nonelementary, then there exists $[T]\in \Lambda_G^{cv}$ such that $[T]\ne [T_\pm(\phi)]$. Then $\lim_{n\to\pm \infty} [T]\phi^n=[T_\pm]\in \Lambda_G^{cv}$. The argument for $[\mu_\pm]$ is the same.
\end{rem}

\begin{defn}\label{defn:delta}
Let $G\le \Out(F_N)$ be a dynamically large subgroup.
Put
\[
\Delta_G^{cv}:=\CVNbar-Z_G^{cv}
\]
and
\[
\Delta_G^{curr}:=\mathbb PCurr(F_N)-Z_G^{curr}
\]
\end{defn}
Note that by construction $\Delta_G^{cv}$ and $\Delta_G^{curr}$ are open $G$-invariant subsets of $\CVNbar$ and of $\mathbb PCurr(F_N)$ accordingly.
\section{Dichotomy}

\begin{conv}\label{conv:G}
Through this section, unless specified otherwise, let $G\le \Out(F_N)$, where $N\ge 3$, be a dynamically large subgroup.
Let $g\in G$ be an atoroidal iwip. Let $[T_+],[T_-]\in \CVNbar$, $[\mu_+],[\mu_-]\in \mathbb PCurr(F_N)$ be the attracting and repelling fixed points of $g$ in $\CVNbar$ and $\mathbb PCurr(F_N)$ accordingly.
For a current $\mu\in Curr(F)$ denote $\langle T_{\pm },\mu\rangle:=\max\{\langle T_-,\mu\rangle, \langle T_+,\mu\rangle\}$.
\end{conv}

\subsection{Basic dichotomy for trees}
\begin{lem}\label{lem:c_n}
Let $T\in cv(F)$ and let $g_n\in \Out(F_N)$ be an infinite sequence of distinct elements. Let $T_\infty\in \cvnbar$, $c_n\ge 0$ be such that $\lim_{n\to\infty} c_nTg_n\to T_\infty$ in $\cvnbar$. Then
$\lim_{n\to\infty}c_n=0$.
\end{lem}
\begin{proof}
Recall that $\CVN=\mathbb P\cvn$ and $\CVNbar=\mathbb P\cvnbar$. Since $\lim_{n\to\infty} [T]g_n=[T_\infty]$ in $\CVNbar$ and since the action of $\Out(F_N)$ on $\CVN$ is properly discontinuous, it follows that $[T_\infty]\in \partial \CVN=\CVNbar-\CVN$. Therefore $F_N$ has nontrivial elements of arbitrarily small translation length with respect to the action on $T_\infty$. On the other hand, the action of $F_N$ on $T$ is free and discrete and therefore there exists $C>0$ such that for every $w\in F_N, w\ne 1$ and for every $n\ge 1$ we have $||w||_{Tg_n}=||g_n(w)||_T\ge C$. The statement of the lemma now follows from point-wise convergence of the translation length functions of $Tg_n$ to that of $T_\infty$.
\end{proof}
\begin{cor}\label{cor:infty}
Let $g_n\in \Out(F_N)$ be an infinite sequence of distinct elements. Then there exists a conjugacy class $[w]$ in $F_N$ such that the set of conjugacy classes $g_n[w]$ is infinite.
\end{cor}
\begin{proof}
Let $T\in \cvn$ be arbitrary.
Choose a limit point $[T_\infty]$ of $[T]g_n$ in $\CVNbar$. Then, after passing to a subsequence, we have $\lim_{n\to\infty} c_nTg_n=T_\infty$. Choose $w\in F_N$ such that $||w||_{T_\infty}>0$. Thus $\lim_{n\to\infty} c_n||g_n(w)||_T=||w||_{T_\infty}>0$. Since by Lemma~\ref{lem:c_n} we have $\lim_{n\to\infty}c_n=0$, it follows that $\lim_{n\to\infty} ||g_n(w)||_T=\infty$, so that the sequence of conjugacy classes $g_n[w]$ contains infinitely many distinct elements.
\end{proof}

\begin{lem}\label{lem:bd}
Let $T\in \cvn$ and $\mu\in Curr(F_N)$. Let $g_n\in \Out(F_N)$ be an infinite sequence of distinct elements. Then one of the following holds:
\begin{enumerate}
\item The sequence $\langle T, g_n\mu\rangle$ is unbounded.
\item For every limit point $[T_\infty]$ of $[T]g_n$ we have $\langle T_\infty,\mu\rangle=0$.
\end{enumerate}
\end{lem}
\begin{proof}
Suppose that (2) fails and there exists a limit point $[T_\infty]$ of $[T]g_n$ such that $\langle T_\infty,\mu\rangle>0$.
After passing to a subsequence $g_{n_i}$ of $g_n$ we have $T_\infty=\lim_{i\to\infty} c_iTg_{n_i}$ for some $c_i\ge 0$. Then
\[
\langle T_\infty,\mu\rangle=\lim_{i\to\infty} \langle c_iTg_{n_i},\mu\rangle = \lim_{i\to\infty} c_i\langle T,g_{n_i}\mu\rangle.
\]
Since by Lemma~\ref{lem:c_n} we have $\lim_{i\to\infty}c_i=0$ and since $\langle T_\infty,\mu\rangle>0$, it follows that $\lim_{i\to\infty} \langle T,g_{n_i}\mu\rangle=\infty$. Therefore the sequence $\langle T, g_n\mu\rangle$ is unbounded, as required.
\end{proof}

Recall that $g\in G$ is an atoroidal iwip and that $[T_\pm]$, $[\mu_\pm]$ are its fixed points in $\CVNbar$ and $\mathbb PCurr(F_N)$.

The following fact is established in ~\cite{KL3}:
\begin{prop}\label{prop:sum}
Let $T_0\in \cvn$ be arbitrary. Then the functions $||.||_{T_0}$
and $||.||_{T_+}+||.||_{T_-}$ on $F_N$ are bi-Lipschitz equivalent.
\end{prop}

\begin{lem}\label{lem:c_n1}
Let $g_n\in \Out(F_N)$ be an infinite sequence of distinct elements. Then there exist a subsequence $g_n'=g_{s_n}$, $s_n\to\infty$, and a tree $T_\ast\in \{T_+,T_-\}$ such that the conclusion of Lemma~\ref{lem:c_n} holds for $T_\ast$ and $g_n'$:

That is, if, after passing to a further subsequence $g_{n_i}'$, we have $\lim_{i\to\infty} c_iT_\ast g_{n_i}'=T_\infty$ in $\cvnbar$ for some $c_i\ge 0$ then $\lim_{i\to\infty}c_i=0$.
\end{lem}
\begin{proof}
By Proposition~\ref{prop:sum}, for any $T_0\in \cvn$ there exists $C>0$ such that
 \[
 ||.||_{T+}+||.||_{T_-}\ge C||.||_{T_0} \text{ on } F_N.
 \]
Corollary~\ref{cor:infty} implies that there exists a nontrivial conjugacy class $[w]$ in $F_N$ such that the set of conjugacy classes $g_n[w]$ is infinite, so that $\lim_{n\to\infty}||g_n(w)||_{T_0}=\infty$. Hence, after passing to a subsequence $g_n'$ of $g_n$,  by the above inequality we conclude that for some $T_\ast\in \{T_+,T_-\}$, we have $\lim_{n\to\infty}||g_n'(w)||_{T_\ast}=\infty$.

Suppose now that $\lim_{i\to\infty} c_iT_\ast g_{n_i}'=T_\infty$.
Since there exists a finite limit
\[
||w||_{T_\infty}=\lim_{i\to\infty} c_i||w||_{T_\ast g_{n_i}'}=\lim_{i\to\infty} c_i||g_{n_i}'(w)||_{T_\ast},
\]
it follows that $\lim_{i\to\infty} c_i=0$, as required.
\end{proof}

\begin{conv}\label{conv:g_n}
For the remainder of this section, unless specified otherwise, we assume that $g_n\in G$ is a (fixed) infinite sequence of distinct elements of a dynamically large subgroup $G\le \Out(F_N)$ and that $g_n'$ and $T_\ast\in \{T_+,T_-\}$ are provided by Lemma~\ref{lem:c_n1}. 
\end{conv}

\begin{cor}\label{cor:unboundedcv}
For any $[\mu] \in
\Delta_G^{curr}$, we have $\lim_{n\to\infty}\langle T_\ast, g_n'\mu\rangle=\infty$. In particular, the sequence $\langle T_\ast, g_n\mu\rangle$ is unbounded.
\end{cor}
\begin{proof}

Suppose that for some $[\mu] \in
\Delta_G^{curr}$, we have $\lim_{n\to\infty}\langle T_\ast, g_n'\mu\rangle\ne \infty$. Then there exists a subsequence $g_n''$ of $g_n'$ such that the sequence
$\langle T_\ast, g_n''\mu\rangle$ is bounded. Since $\CVNbar$ is compact, after passing to a further subsequence $g_{n_i}$ of $g_n''$, there exist $c_i\ge 0$ and $T_\infty\in \cvnbar$ such that $\lim_{i\to\infty} c_iT_\ast g_{n_i}=T_\infty$. By Lemma~\ref{lem:c_n1} we have $\lim_{i\to\infty} c_i=0$. Therefore, as in the proof of Lemma~\ref{lem:bd}, either the sequence $\langle T_\ast,g_{n_i}\mu\rangle$ is unbounded or $\langle T_\infty,\mu\rangle=0$. The former is impossible since by assumption the sequence
$\langle T_\ast, g_n''\mu\rangle$ is bounded. Thus $\langle T_\infty,\mu\rangle=0$.
However, $[T_{\pm }]\in \Lambda_G^{cv}$ by Remark~\ref{rem:poles} and hence $[T_{\infty}]\in \Lambda_G^{cv}$. Therefore the possibility that $\langle T_\infty,\mu\rangle=0$ is ruled out by the assumption that $[\mu]\in \Delta_G^{curr}$, yielding a contradiction.
\end{proof}

\subsection{Basic dichotomy for currents}

We now want to dualize the previous arguments for the case of geodesic currents. The difficulty here is that currents analogues of Lemma~\ref{lem:c_n} and Lemma~\ref{lem:c_n1} are not readily available. We get around this problem by using the results of the previous steps to obtain such analogues first.

\begin{lem}\label{lem:currscalars}
Let $[\mu]\in \Delta_G^{curr}$.  Suppose that for a subsequence $g_{n_i}$ of $g_n'$ there are $c_n\ge 0$ and $\mu_\infty\in Curr(F_N)$ be such that $\lim_{i\to\infty} c_ig_{n_i}\mu=\mu_\infty$. Then $\lim_{i\to\infty} c_i=0$.
\end{lem}
\begin{proof}
We have
\[
\langle T_\ast,\mu_\infty\rangle=\lim_{i\to\infty} \langle T_\ast, c_ig_{n_i}\mu\rangle=\lim_{i\to\infty} c_i \langle T_\ast,g_{n_i}\mu\rangle.
\]
 By Corollary~\ref{cor:unboundedcv} we have $\lim_{i\to\infty}\langle T_\ast, g_{n_i}\mu\rangle=\infty$.
Hence, from the above equation, we conclude that $\lim_{i\to\infty} c_i=0$, as required.
\end{proof}

We now obtain an analogue of Lemma~\ref{lem:bd} for currents, with the important difference however, that, as specified in Convention~\ref{conv:g_n}, here we do need to assume that $g_n\in G$ while the proof of Lemma~\ref{lem:bd} works for an arbitrary infinite sequence $g_n\in \Out(F_N)$. Similarly, we had to assume that $g_n\in G$ in Lemma~\ref{lem:currscalars} while this assumption was not needed in Lemma~\ref{lem:c_n}.
Recall that $g_n\in G$, as well as $g_n'$ and $T_\ast$ are as in Convention~\ref{conv:g_n}.

\begin{lem}\label{lem:bd-curr}
 Let $T\in \cvnbar$ be arbitrary and let $[\mu]\in \Delta_G^{curr}$. Then one of the following holds:
\begin{enumerate}
\item The sequence $\langle Tg_n', \mu\rangle$ is unbounded.
\item We have $\langle T,\mu_\infty\rangle=0$ for any limit point $[\mu_\infty]$ of the sequence $[g_n'\mu]$ in $\mathbb PCurr(F_N)$.
\end{enumerate}
\end{lem}
\begin{proof}
Suppose that (2) fails and there exist $c_i\ge 0$ and $\mu_\infty\in Curr(F_N)$ be such that, after possibly passing to a subsequence $g_{n_i}$ of $g_n'$, we have $\lim_{i\to\infty} c_ig_{n_i}\mu=\mu_\infty$ and such that $\langle T,\mu_\infty\rangle>0$.
Then
\[
\langle T,\mu_\infty\rangle=\lim_{i\to\infty} \langle T, c_ig_{n_i}\mu\rangle=\lim_{i\to\infty} c_i\langle Tg_{n_i}, \mu\rangle.
\]
Since by Lemma~\ref{lem:currscalars} we have $\lim_{i\to\infty}c_i=0$, and since by assumption $\langle T,\mu_\infty\rangle>0$, it follows that the sequence $\langle Tg_{n}', \mu\rangle$ is unbounded.
\end{proof}

\begin{lem}\label{lem:currpmscalars}
There exists $\mu\in \{\mu_\pm\}$ with the following property. Whenever, after passing to a further subsequence $g_{n_i}$ of $g_n'$, we have $c_{i}g_{n_i}\mu\to\mu_\infty\ne 0$ for some $c_{i}\ge 0$, then $\lim_{i\to\infty} c_{i}=0$.
\end{lem}

\begin{proof}

Suppose there exist $c_i, g_{n_i'}$ and $\mu_\infty$ as above such that $c_{i}g_{n_i}\mu\to\mu_\infty\ne 0$ but $\lim_{i\to\infty} c_{i}\ne 0$.
Then after passing to a further subsequence, we may assume that $c_{n_i}\ge c>0$ for all $i\ge 1$.

Let $T\in \cvn$.
First, we claim that either the sequence $\langle T, g_{n_i}\mu_+\rangle$ or the sequence $\langle T,g_{n_i}\mu_-\rangle$ is unbounded. Suppose not, that is, that they are both  bounded. Choose a limit point $[T_\infty]$ of $[Tg_{n_i}]$. Then by Lemma~\ref{lem:bd} we have $\langle T_\infty,\mu_+\rangle=\langle T_\infty, \mu_-\rangle=0$. However, this is impossible by Proposition~\ref{prop:zero}. Thus there is $\mu\in \{\mu\pm\}$ such that the sequence $\langle T, g_{n_i}\mu\rangle$ is unbounded. After taking a further subsequence of $g_{n_i}$, we may thus assume that $\lim_{i\to\infty} \langle T, g_{n_i}\mu\rangle=\infty$. Recall that $c_i\ge c>0$ for every $i\ge 1$.
Now,
\[
\langle T,\mu_\infty\rangle=\lim_{i\to\infty} \langle T, c_{i}g_{n_i}\mu\rangle=\lim_{i\to\infty} c_i\langle T, g_{n_i}\mu\rangle
\]
which yields a contradiction since  $\lim_{i\to\infty} \langle T, g_{n_i}\mu\rangle=\infty$ and $c_i\ge c>0$ for every $i\ge 1$.
\end{proof}

\begin{cor}\label{cor:unboundedcurrents}
Let $\mu\in \{\mu_\pm\}$ be provided by Lemma~\ref{lem:currpmscalars}. Let $[T]\in \Delta_G^{cv}$.
Then the sequence $\langle Tg_n,\mu\rangle$ is unbounded.
\end{cor}
\begin{proof}
We first pass to the subsequence $g_n'$ of $g_n$ provided by Lemma~\ref{lem:c_n1} so that $\mu$ has the properties guaranteed by Lemma~\ref{lem:currpmscalars}. Let $[\mu_\infty]$ be a limit point of $[g_n'\mu]$. Thus, for a subsequence $g_{n_i}$ of $g_n'$ and for some $c_i\ge 0$ we have $\lim_{n\to\infty} c_ig_{n_i}\mu=\mu_\infty$. By Lemma~\ref{lem:currpmscalars} we have $\lim_{i\to\infty} c_i=0$. Note also that $\mu\in \{\mu_\pm\}$ which implies that $[\mu]\in \Lambda_G^{curr}$ and hence $[\mu_\infty]\in \Lambda_G^{curr}$.

Now let $[T]\in \Delta_G^{cv}$. Then
\[
\langle T,\mu_\infty\rangle =\lim_{i\to\infty} c_i \langle T,g_{n_i}\mu\rangle=\lim_{i\to\infty} c_i \langle Tg_{n_i},\mu\rangle.
\]
Note that $\langle T,\mu_\infty\rangle>0$ since $[\mu_\infty]\in \Lambda_G^{curr}$ and $[T]\in \Delta_G^{cv}$.
By Lemma~\ref{lem:currpmscalars} we have $\lim_{i\to\infty} c_i=0$ which now implies that $\lim_{i\to\infty} \langle Tg_{n_i},\mu\rangle=\infty$. Thus the sequence $\langle Tg_n,\mu\rangle$ is unbounded, as claimed.
\end{proof}

\begin{cor}\label{cor:c_n2}
Let $[T]\in \Delta_G^{cv}$ and let $c_i\ge 0$, and let $g_{n_i}$ be a subsequence of $g_n'$ such that $\lim_{i\to\infty} c_i Tg_{n_i}=T_\infty\in \cvnbar$. Then $\lim_{i\to\infty} c_i=0$.
\end{cor}
\begin{proof}
Let $\mu\in \{\mu_\pm\}$ have the properties guaranteed by Lemma~\ref{lem:currpmscalars}.
Suppose that $\lim_{i\to\infty} c_i\ne 0$. Then, after passing to a further subsequence, we have $c_i\ge c>0$ for all $i\ge 1$.
After passing to a further subsequence, we may assume that $c_i'g_{n_i}\mu\to\mu_{\infty}$. We have $\lim_{i\to\infty} c_i'=0$. The proof of Corollary~\ref{cor:unboundedcurrents} shows that in this case $\lim_{i\to\infty} \langle Tg_{n_i},\mu\rangle=\infty$ since
\[
\langle T,\mu_\infty\rangle =\lim_{i\to\infty} c_i' \langle T,g_{n_i}\mu\rangle=\lim_{i\to\infty} c_i' \langle Tg_{n_i},\mu\rangle.
\]

Then
\[
\langle T_\infty, \mu\rangle=\lim_{i\to\infty} c_i \langle Tg_{n_i},\mu\rangle=\infty,
\]
since $c_i\ge c>0$ and since $\lim_{i\to\infty} \langle Tg_{n_i},\mu\rangle=\infty$. This yield a contradiction since the value of intersection form is always finite.
\end{proof}

\begin{prop}Let $G\le \Out(F_N)$ be a dynamically large subgroup.
Then the following hold:
\begin{enumerate}
\item Let $[T]\in \Delta_G^{cv}$ and let $[T_\infty]$ be a limit point of $G[T]$. Then $[T_\infty]\in Z_G^{cv}$.
\item Let $[\mu]\in \Delta_G^{curr}$ and let $[\mu_\infty]$ be a limit point of $G[\mu]$. Then $[\mu_\infty]\in Z_G^{curr}$.
\end{enumerate}
\end{prop}
\begin{proof}
Let $[T]\in \Delta_G^{cv}$ and let $[T_\infty]$ be a limit point of $G[T]$. Then there is an infinite sequence of distinct elements $g_n\in G$ such that $\lim_{n\to\infty} [T]g_n=[T_\infty]$. Let $g_n'$ be the subsequence of $g_n$ provided by Lemma~\ref{lem:c_n1}. Then for some subsequence $g_{n_i}$ of $g_n'$ and for some $c_i\ge 0$ we have $\lim_{i\to\infty} c_iTg_{n_i}=T_\infty$. By Corollary~\ref{cor:c_n2} we have $\lim_{i\to\infty} c_i=0$. We now apply Lemma~\ref{lem:currpmscalars} to the sequence $h_i=g_{n_i}^{-1}$. Lemma~\ref{lem:currpmscalars} implies that after passing to a further subsequence for some $\mu\in \{\mu_\pm\}$ we have $\lim_{i\to\infty} c_i'g_{n_i}^{-1}\mu=\mu_\infty\ne 0$ and $\lim_{i\to\infty}c_i'=0$.

Therefore
\[
\langle T_\infty,\mu_\infty\rangle=\lim_{i\to\infty} c_ic_i'\langle Tg_{n_i}, g_{n_i}^{-1}\mu\rangle=\lim_{i\to\infty} c_ic_i'\langle T, \mu\rangle=0.
\]
By construction we have $[\mu]=[\mu_\pm]\in \Lambda_G^{curr}$ and hence $[\mu_\infty]=\lim_{i\to\infty} g_{n_i}^{-1}[\mu]\in \Lambda_G^{curr}$. Therefore by definition of the zero-set we have $[T_\infty]\in Z_G^{cv}$, as required.

The proof of part (2) of the proposition is similar. The use of Corollary~\ref{cor:c_n2} and Lemma~\ref{lem:currpmscalars} has to be replaced by the use of Lemma~\ref{lem:currscalars} and Lemma~\ref{lem:c_n1}. We leave the details to the reader.
\end{proof}
\section{Domains of discontinuity}

\begin{conv}
For the remainder of this section let $G\le \Out(F_N)$ be a dynamically large subgroup, let $\phi\in G$ be an atoroidal iwip and let $[T_\pm]$, $[\mu_\pm]$ be the fixed points of $\phi$ in $\CVNbar$ and in $\mathbb PCurr(F_N)$ accordingly.
Recall that for $T\in \cvnbar$ we denote $\langle T, \mu_\pm\rangle=\langle T,\mu_+\rangle+\langle T,\mu_-\rangle$. Similarly, for $\mu\in Curr(F_N)$ we denote $\langle T_\pm, \mu\rangle=\langle T_+,\mu\rangle+\langle T_-,\mu\rangle$.
\end{conv}

\begin{notation}
Denote
\[
D_G^{curr}=\{[\mu]\in \Delta_G^{curr}: \langle T_\pm, \mu\rangle \le \langle T_\pm, g\mu\rangle \text{ for every } g\in G\}
\]
and
\[
D_G^{cv}=\{[T]\in \Delta_G^{cv}: \langle T, \mu_\pm\rangle \le \langle gT, \mu_\pm\rangle \text{ for every } g\in G\}
\]
\end{notation}

\begin{lem}\label{lem:fd}
We have $G D_G^{curr}=\Delta_G^{curr}$ and $G D_G^{cv}=\Delta_G^{cv}$.
\end{lem}

\begin{proof}
Let $[\mu]\in \Delta_G^{curr}$. Then for any $C\ge 1$ the set
\[
\{g\in G: \langle g\mu, T_\pm\rangle\le C\}
\]
is finite by Corollary~\ref{cor:unboundedcv}. Hence there exists $g_0\in G$ with $\langle g_0\mu, T_\pm\rangle=\min_{g\in G} \langle g\mu, T_\pm\rangle$, that is, $g_0[\mu]\in D_G^{curr}$. Therefore $G D_G^{curr}=\Delta_G^{curr}$, as required.
The proof that $G D_G^{cv}=\Delta_G^{cv}$ is similar, with the use of Corollary~\ref{cor:unboundedcv} replaced by Corollary~\ref{cor:unboundedcurrents}.
\end{proof}

\begin{lem}\label{lem:pd}
For any compact set $K'\subseteq\Delta_G^{curr}$ and for any compact set $K\subseteq \Delta_G^{cv}$ we have
\[
\{g\in G: K'\cap gD_G^{curr}\ne \emptyset\} \text{ is finite}
\]
and
\[
\{g\in G: K\cap gD_G^{cv}\ne \emptyset\} \text{ is finite.}
\]
\end{lem}
\begin{proof}
Suppose that there exists an infinite sequence of distinct elements $g_n\in G$ such that $K\cap g_nD_G^{cv}\ne \emptyset$ for all $n\ge 1$.
Then there is also a sequence $[T_n]\in K\cap g_nD_G^{cv}$. After passing to a subsequence, we may assume that $[T_n]\to[T_\infty]$ in $\CVNbar$. Moreover, after choosing the projective representatives of $[T_n]$ appropriately, we may assume that $T_n\to T_\infty$ as $n\to\infty$. Note that $[T_\infty]\in K$ since $K$ is compact.

By Lemma~\ref{lem:currpmscalars}, after passing to a further subsequence, there exist $\mu\in \{\mu_\pm\}$, $c_n\ge 0$ and $0\ne\mu_\infty\in Curr(F_N)$ such that $\lim_{n\to\infty} c_ng_n\mu=\mu_\infty$ and such that $\lim_{n\to\infty} c_n=0$.
We have
\begin{gather*}
\langle T_\infty,\mu_\infty\rangle =\lim_{n\to\infty} c_n \langle T_n, g_n\mu\rangle\le \lim_{n\to\infty} c_n \langle T_n, g_n\mu_\pm\rangle =\\
=\lim_{n\to\infty} c_n \langle g_n^{-1}T_n,\mu_\pm\rangle \le \qquad \text{since } [T_n]\in g_nD_G^{cv} \\
\le \lim_{n\to\infty} c_n \langle T_n,\mu_\pm\rangle=\lim_{n\to\infty} c_n \lim_{n\to\infty} \langle T_n,\mu_\pm\rangle=0\cdot \langle T_\infty, \mu_\pm\rangle=0.
\end{gather*}
Thus $\langle T_\infty,\mu_\infty\rangle=0$. Recall that $[T_\infty]\in K\subseteq \Delta_G^{cv}$. Also, $[\mu]\in\{[\mu_\pm]\}$ and hence $[\mu]\in \Lambda_G^{curr}$ and therefore $[\mu_\infty]\in \Lambda_G^{curr}$. The fact that $\langle T_\infty,\mu_\infty\rangle=0$ now gives a contradiction with the definition of $\Delta_G^{cv}$.

The proof that $\{g\in G: K'\cap gD_G^{curr}\ne \emptyset\}$ is finite is similar. The only difference is that in the argument the use of  Lemma~\ref{lem:currpmscalars} has to be replaced by the use of Lemma~\ref{lem:c_n1}.
\end{proof}

\begin{thm}\label{thm:pd}
Let $G\le \Out(F_N)$ be a dynamically large subgroup. Then the actions of $G$ on $\Delta_G^{cv}$ and $\Delta_G^{curr}$ are properly discontinuous.
\end{thm}
\begin{proof}
We will show that the action of $G$ on  $\Delta_G^{curr}$ is properly discontinuous. The argument for $\Delta_G^{cv}$ is similar.

Let $K'\subseteq \Delta_G^{curr}$ be a compact subset and suppose that $\{g\in G: K'\cap gK'\ne\emptyset\}$ is infinite. Let $g_n\in G$ be an infinite sequence of distinct elements such that $K'\cap g_nK'\ne\emptyset$ for every $n\ge 1$.

By Lemma~\ref{lem:pd} there exists a finite collection $h_1,\dots, h_t\in G$ such that
\[
\{h_1,\dots, h_t\}=\{g\in G: K'\cap gD_G^{curr}\ne \emptyset\}.
\]
Moreover, since by Lemma~\ref{lem:fd} $GD_G^{curr}=\Delta_G^{curr}$, it follows that $K'\subseteq \cup_{i=1}^t h_i D_G^{curr}$.

For every $n\ge 1$ there exists $[\mu_n]\in K'\cap g_nK'$ and hence $[\mu_n]\in K'\cap \cup_{i=1}^t g_nh_i D_G^{curr}$.
Therefore for every $n\ge 1$ there exists $i_n, 1\le i_n\le t$ such that $K'\cap g_nh_{i_n}D_G^{curr}\ne \emptyset$. It follows that
\[
\{g_nh_{i_n}: n\ge 1\}\subseteq \{h_1,\dots, h_t\}
\]
which yields a contradiction since by assumption the sequence $g_n, n\ge 1$ consists of infinitely many distinct elements of $G$.
\end{proof}

\begin{rem}\label{rem:maximal}
The domains of discontinuity constructed in Theorem~\ref{thm:pd} are
not necessarily maximal. Thus for the case of $G=\Out(F_N)$
Guirardel~\cite{Gui1} constructed an open $G$-invariant subset $\mathcal
O_N\subset \CVNbar$ such that $\Out(F_N)$ acts on $\mathcal
O_N$ properly discontinuously. In Guirardel's construction,
$\CVN\subsetneq \mathcal O_N$. One can see that for every $[T]\in
\mathcal O_N-\CVN$ we have $[T]\in Z_{Out(F_N)}^{cv}$ and hence
$[T]\not\in \Delta_G^{cv}$. The reason for this is that every such $T$
corresponds to a (non-free) simplicial action of $F_N$ with trivial
edge stabilizers (plus some additional conditions on the quotient
graph that are not relevant here) and therefore there exists a
primitive element $a$ of $F_N$ such that $a$ fixes a vertex in
$T$. For every such $a$ we have $[\eta_a]\in \Lambda_{Out(F_N)}^{curr}$
(see \cite{KL1}) and $\langle T,\eta_a\rangle=||a||_T=0$, so that
indeed $[T]\in Z_{Out(F_N)}^{cv}$, as claimed.
\end{rem}

Let $Curr_+(F_N)$ be the set of all $\mu\in Curr(F_N)$ with full support and let $\mathbb PCurr_+(F_N)=\{[\mu]: \mu\in Curr_+(F_N)\}$.
Note that $\mathbb PCurr_+(F_N)$ is an open $\Out(F_N)$-invariant subset of $\mathbb PCurr(F_N)$.

\begin{cor}\label{cor:+}
Let $N\ge 3$. Then the action of $\Out(F_N)$ on $\mathbb PCurr_+(F_N)$ is properly discontinuous.
\end{cor}
\begin{proof}
Note that for $N\ge 3$ the group $G=\Out(F_N)$ is dynamically large.
By the main result of~\cite{KL3} we have $\langle T,\mu\rangle>0$ for any $T\in \cvnbar$ and any $\mu\in Curr_+(F_N)$. Therefore, by Definition~\ref{defn:delta}, we have $\mathbb PCurr_+(F_N)\subseteq \Delta_{Out(F_N)}^{curr}$. Hence the action of $\Out(F_N)$ on $\Delta_{Out(F_N)}^{curr}$ is properly discontinuous by Theorem~\ref{thm:pd} and therefore the action of $\Out(F_N)$ on $\mathbb PCurr_+(F_N)$ is properly discontinuous as well.
\end{proof}

Note, however, that in the proof of Corollary~\ref{cor:+} the containment $\mathbb PCurr_+(F_N)\subsetneq \Delta_{Out(F_N)}^{curr}$ is proper. Recall that a current $\mu\in Curr(F_N)$ is called \emph{filling} if for every $T\in \cvnbar$ we have $\langle T,\mu\rangle >0$. The same argument as in the proof of Corollary~\ref{cor:+} shows that if $\mu\in Curr(F_N)$ is filling then $[\mu]\in \Delta_{Out(F_N)}^{curr}$. It was proved in~\cite{KL3} that there exist filling rational currents $\eta_g\in Curr(F_N)$ and in fact, the property of being filling for a rational current is, in a sense, ``generic". However, a rational current never has full support, so if $\eta_g$ is filling, then $\eta_g\in \Delta_{Out(F_N)}^{curr}-\mathbb PCurr_+(F_N)$.

\end{document}